\documentclass[11pt,reqno]{amsart}
\usepackage{amsmath}
\usepackage{amssymb}
\usepackage{enumerate}

\newtheorem{theorem}{Theorem}[section]
\newtheorem{lemma}[theorem]{Lemma}
\newtheorem{corollary}[theorem]{Corollary}

\theoremstyle{definition}
\newtheorem{remark}[theorem]{Remark}

\newtheorem{question}{Question}
\newtheorem{example}[theorem]{Example}

\numberwithin{equation}{section}

\allowdisplaybreaks[1]

\begin{document}

\title[Toeplitz operators on polyanalytic functions]{Toeplitz operators on Bergman spaces of polyanalytic functions}

\author{{\v Z}eljko {\v C}u{\v c}kovi{\'c}}
\address{Department of Mathematics, Mail Stop 942, University of Toledo, Toledo, OH 43606}
\email{zcuckovi@math.utoledo.edu}

\author{Trieu Le}
\address{Department of Mathematics, Mail Stop 942, University of Toledo, Toledo, OH 43606}
\email{trieu.le2@utoledo.edu}


\begin{abstract} We study algebraic properties of Toeplitz operators on Berg-man spaces of polyanalytic functions on the unit disk. We obtain results on finite-rank commutators and semi-commutators of Toeplitz operators with harmonic symbols. We also raise and discuss some open questions.
\end{abstract}
\maketitle


\section{Introduction}\label{S:intro}
The last twenty years have brought a lot of progress in the study of spaces of analytic functions and their operators. Particularly interesting to us are Bergman spaces and Toeplitz operators acting on them. One can notice the development of the field in several directions. One promising direction is to study these operators on Bergman spaces of different domains in $\mathbb{C}^n$. The other generalization concerns harmonic Bergman spaces. We choose to go to a different generalization, namely we would like to study algebraic properties of Toeplitz operators on Bergman spaces of polyanalytic functions on the unit disk $\mathbb{D}$ in the complex plane. For $n\geq 1$, a function $h$ is called $n$-analytic (or polyanalytic of order $n$) on $\mathbb{D}$ if it satisfies the generalized Cauchy-Riemann equation $\frac{\partial^n h}{\partial\bar{z}^n}=0$. It is clear that $h$ is $n$-analytic if and only if there are analytic functions $h_0, \ldots, h_{n-1}$ on $\mathbb{D}$ such that $$h(z)=h_0(z)+h_1(z)\bar{z}+\cdots+h_{n-1}(z)\bar{z}^{n-1},$$ for all $z\in\mathbb{D}$. 

Throughout the paper, we use $L^p(\mathbb D)$ ($1\leq p\leq\infty$) to denote the space $L^p(\mathbb{D},\mathrm{dA})$, where $\mathrm{dA}$ is the normalized Lebesgue area measure on $\mathbb D$. The $n$-analytic Bergman space $A^2_{n}=A^2_{n}(\mathbb{D})$ consists of all $n$-analytic functions which also belong to the Hilbert space $L^2(\mathbb{D})$ (with inner product denoted by $\langle\cdot,\cdot\rangle$). In other words, a function $h$ belongs to $A^2_n$ if and only if $h$ is $n$-analytic and $$\|h\|=\langle h,h\rangle^{1/2} = \Big(\int_{\mathbb D}|h(z)|^2\mathrm{dA}(z)\Big)^{1/2}<\infty.$$
The space $A_1^2$ is the familiar Bergman space of the unit disk. It is well known that $A_1^2$ is a closed subspace of $L^2(\mathbb{D})$. It turns out that the same thing holds for all $n\geq 1$. 

One of the first studies of $A^2_n$ appeared in Ko{\v{s}}elev's paper \cite{Koshelev1977}. We also refer the interested reader to Ramazanov's papers \cite{RamazanovMZ2002,RamazanovMZ2004} and Balk's book \cite{Balk1991} and the references therein for more details about the structure of these spaces. 

For any function $f\in L^1(\mathbb{D},\mathrm{dA})$, the Toeplitz operator $T_f$ is an integral operator defined by the formula
\begin{equation*}
(T_fh)(z) = \int_{\mathbb D}f(w)h(w)\bar{K}_z(w)\mathrm{dA}(w),
\end{equation*}
for $h\in A^2_n$ for which the integral on the right hand side is well defined for all $z\in\mathbb{D}$. Here $K_z(w)$ is the reproducing kernel function of $A^2_n$. When $f$ is a bounded function, $T_fh$ is defined for all $h\in A^2_n$ and it is also given by $T_fh=P(fh)$, where $P$ is the orthogonal projection from $L^2(\mathbb{D})$ onto $A^2_n$.

The first study of Toeplitz operators on $A^2_n$ seems to be Wolf's paper \cite{WolfMN1994}. In the paper, Wolf considered Toeplitz operators $T_f$ where the function $f$ is continuous on the closed unit disk. He showed that such an operator is unitarily equivalent to a compact perturbation of the direct sum of $n$ copies of the Hardy-Toeplitz operator whose symbol is the restriction of $f$ on the unit circle. As a consequence, a criterion for Fredholmness and a Fredholm index formula for $T_f$ were obtained.

The goal of the current paper is to use the reproducing kernels and Berezin transform to understand the algebraic and quantitative properties of Toeplitz operators on $A^2_n$. We show that certain analytic Bergman space results remain true in the new setting. In particular, we generalize the result of the first author \cite{CuckovicIEOT2007} on finite rank perturbations of products of Toeplitz operators with harmonic symbols. Such a generalization is non trivial since the new Berezin transform is no longer one-to-one when $n\geq 2$, unlike the analytic Bergman space case when $n=1$. We hope to be able to deepen our understanding of operators on these function spaces in our future work.

\section{Preliminaries}\label{S:prel}

In this section we discuss in more detail the spaces $A^2_n$ and their reproducing kernel functions.

We begin by showing that the evaluation functional at each point in the unit disk is bounded on $A^2_n$. Recall that the pseudo-hyperbolic distance between two points $z$ and $w$ in the disk is defined by $|\varphi_z(w)|$, where $\varphi_{z}(w)=\frac{z-w}{1-w\bar{z}}$ is the Moebius transformation of the disk that interchanges $0$ and $z$. For $z\in\mathbb{D}$ and $0<r<1$, let $E(z,r)=\{w: |\varphi_z(w)|<r\}$ be the pseudo-hyperbolic disk of radius $r$ centered at $z$. It is well known (see \cite[Proposition 4.4]{ZhuAMS2007}) that $E(z,r)$ is the Euclidean disk of radius $R$ centered at $c$, where
$$c=\dfrac{1-r^2}{1-r^2|z|^2}z,\quad \text{ and }\quad R=\dfrac{1-|z|^2}{1-r^2|z|^2}r.$$
Since $\tilde{R}=R-|c-z|=\dfrac{(1-|z|^2)r}{1+r|z|}>0$, we see that $E(z,r)$ contains the Euclidean disk $\{w: |w-z|<\tilde{R}\}$.

In \cite[Lemma 5]{PavlovicJMAA1997}, Pavlovi{\'c} showed that there is a constant $C>0$ such that for any $h\in A^2_n$, $z\in\mathbb{D}$ and $0<\delta<1-|z|$, 
\begin{align}\label{Eqn:meanValues}
|h(z)|^2 & \leq \dfrac{C}{\delta^2}\int_{|w-z|<\delta}|h(w)|^2\mathrm{dA}(w).
\end{align}
Pavlovi{\'c} in fact studied spaces of polyharmonic functions but since any $n$-analytic functions is also $n$-harmonic, his lemma may be used in the current setting. Applying \eqref{Eqn:meanValues} with $\delta=\tilde{R}$ and using the inclusion $\{w: |w-z|<\tilde{R}\}\subset E(z,r)$, we obtain
\begin{align}\label{Eqn:meanBergman}
|h(z)|^2 & \leq\dfrac{C(1+r|z|)^2}{r^2(1-|z|^2)^2}\int_{|w-z|<\tilde{R}}|h(w)|^2\mathrm{dA}(w)\notag\\
& \leq\dfrac{C(1+r)^2}{r^2(1-|z|^2)^2}\int_{E(z,r)}|h(w)|^2\mathrm{dA}(w)\\
& = \dfrac{C_r}{(1-|z|^2)^2}\int_{E(z,r)}|h(w)|^2\mathrm{dA}(w),\notag
\end{align}
where $C_r$ is independent of $h$ and $z$. This implies that $A^2_n$ is a closed subspace of $L^2(\mathbb{D})$ and that the functional $h\mapsto h(z)$ is bounded on $A^2_n$ for any $z\in\mathbb{D}$. Therefore there is a function $K_z$ in $A^2_n$ such that $h(z)=\langle h,K_z\rangle$ for $z\in\mathbb{D}$. The function $K(w,z)=K_z(w)$ for $(z,w)\in\mathbb{D}\times\mathbb{D}$ is called, as usual, the kernel function for $A^2_n$. It can be showed (see \cite{Koshelev1977}) that 
\begin{align*}
K_{z}(w) & = \dfrac{n}{(1-w\bar{z})^{2n}}\sum_{j=0}^{n-1}(-1)^{j}\binom{n}{j+1}\binom{n+j}{n}|1-w\bar{z}|^{2(n-1-j)}|w-z|^{2j}\\
& = \dfrac{n|1-w\bar{z}|^{2(n-1)}}{(1-w\bar{z})^{2n}}\sum_{j=0}^{n-1}(-1)^{j}\binom{n}{j+1}\binom{n+j}{n}|\varphi_{z}(w)|^{2j}.
\end{align*}
Using properties of reproducing kernel functions we obtain, for $z\in\mathbb{D}$,
$\|K_z\|^2 = \langle K_z,K_z\rangle = K_z(z)=n^2/(1-|z|^2)^2.$ The normalized reproducing kernel at $z$ is defined by
\begin{align}\label{Eqn:normalizedKernel}
k_z(w) & = \frac{K_z(w)}{\|K_z\|} = \frac{1-|z|^2}{n}K_z(w)\\
& =  \dfrac{(1-|z|^2)|1-w\bar{z}|^{2(n-1)}}{(1-w\bar{z})^{2n}}\sum_{j=0}^{n-1}(-1)^{j}\binom{n}{j+1}\binom{n+j}{n}|\varphi_{z}(w)|^{2j}.\notag
\end{align}
Using the formula $1-|\varphi_z(w)|^2=\frac{(1-|z|^2)(1-|w|^2)}{|1-w\bar{z}|^2}$ we obtain
{\allowdisplaybreaks
\begin{align*}
|k_z(w)| & = \dfrac{1-|z|^2}{|1-w\bar{z}|^2}\Big|\sum_{j=0}^{n-1}(-1)^{j}\binom{n}{j+1}\binom{n+j}{n}|\varphi_{z}(w)|^{2j}\Big|\\
& =\frac{1-|\varphi_z(w)|^2}{1-|w|^2}\Big|\sum_{j=0}^{n-1}(-1)^{j}\binom{n}{j+1}\binom{n+j}{n}|\varphi_{z}(w)|^{2j}\Big|\\
& = \frac{1-|\varphi_z(w)|^2}{1-|w|^2}\Big|n+\sum_{j=1}^{n-1}(-1)^{j}\binom{n}{j+1}\binom{n+j}{n}|\varphi_{z}(w)|^{2j}\Big|.
\end{align*}}%
The last formula shows that $(1-|w|^2)|k_z(w)|\to n$ as $|\varphi_z(w)|\to 0$, that is, for any $\epsilon>0$, there is a $\delta>0$ such that $\big|(1-|w|^2)|k_z(w)|-n\big|<\epsilon$ for all $w,z\in\mathbb{D}$ satisfying $|\varphi_z(w)|<\delta$. This implies that there is a number $0<r<1$ such that $(1-|w|^2)|k_z(w)|>n/2\geq 1/2$ whenever $|\varphi_z(w)|<r$. Using this we obtain the following estimate.
\begin{lemma}\label{L:integralEstimate}
There is a constant $c>0$ such that for any function $h\in A^2_n$ and any $w\in\mathbb{D}$,
\begin{align*}
|h(w)|^2\leq c\int_{\mathbb{D}}|h(z)|^2|k_z(w)|^2\mathrm{d}A(z).
\end{align*}
\end{lemma}
\begin{proof}
Let $r$ be the number mentioned in the paragraph preceding the lemma. Since $|k_z(w)|>\frac{1}{2(1-|w|^2)}$ for all $z\in E(w,r)$, using \eqref{Eqn:meanBergman}, we have
\begin{align*}
\int_{\mathbb{D}}|h(z)|^2|k_z(w)|^2\mathrm{d}A(z) & \geq\frac{1}{4(1-|w|^2)^2}\int_{E(w,r)}|h(z)|^2\mathrm{d}A(z)\geq \frac{|h(w)|^2}{4C_r}.
\end{align*}
Set $c=4C_r$. The conclusion of the lemma follows.
\end{proof}

The last lemma of this section, which is well known in the case of analytic Bergman space, is a consequence of the fact that $\|K_z\|\to\infty$ as $|z|\uparrow 1$ and $n$-analytic polynomials are dense in $A^2_n$.
\begin{lemma}\label{L:behaviorPoly} We have $k_z\longrightarrow 0$ weakly as $|z|\uparrow 1$. As a result, for any $f\in A^2_n$, $$\lim\limits_{|z|\uparrow 1}(1-|z|^2)f(z)=n\cdot\lim\limits_{|z|\uparrow 1}\langle f, k_z\rangle = 0.$$
\end{lemma}


\section{Weighted Berezin transform}\label{S:wBTransform}

For any real number $\alpha\geq 0$, the weighted Berezin transform $B_{\alpha}$ is defined by
\begin{align*}
B_{\alpha}u(z) & = (\alpha+1)\int_{\mathbb{D}}u(\varphi_{z}(w))(1-|w|^2)^{\alpha}\mathrm{dA}(w)
\end{align*}
for $u$ in $L^1(\mathbb{D})$ and $z$ in $\mathbb{D}$. Note that $B_0$ is just the standard unweighted Berezin transform.

It follows from the mean value theorem for harmonic functions that if $u\in L^1(\mathbb{D})$ is harmonic, then $B_{\alpha}u=u$. Fixed points of $B_{0}$ (even in the setting of the unit ball in higher dimensions) were studied by Ahern, Flores and Rudin in \cite{AhernJFA1993}. In the current setting of the unit disk, their results say that all fixed points of $B_0$ are integrable harmonic functions. This is not true for all $B_{\alpha}$ and the interested reader is referred to \cite[Chapter 6]{ZhuAMS2007} and the references therein.

By a change of variable, we see that
\begin{align*}
B_{\alpha}u(z) & = (\alpha+1)\int_{\mathbb{D}}u(w)\dfrac{(1-|z|^2)^{2+\alpha}(1-|w|^2)^{\alpha}}{|1-z\bar{w}|^{2\alpha+4}}{\rm d}A(w).
\end{align*}
This shows that $B_{\alpha}u$ is a real analytic function on $\mathbb{D}$, hence it is infinitely differentiable.

There is a relation between $B_{\alpha+1}$ and $B_{\alpha}$ that involves the invariant Laplacian $\tilde{\Delta}$. Recall that $\Delta = \frac{\partial^2}{\partial x^2}+\frac{\partial^2}{\partial y^2}=4\frac{\partial^2}{\partial z\partial\bar{z}}$ is the ordinary Laplacian. The invariant Laplacian is defined by $(\tilde{\Delta}u)(z)=(1-|z|^2)^2(\Delta u)(z)$ for any twice differentiable function $u$ on $\mathbb{D}$. 

\begin{lemma}[{\cite[Proposition 2.4]{AhernJFA1993}} and also {\cite[Lemma 6.23]{ZhuAMS2007}}]\label{L:weightedBerezin} If $u\in L^1(\mathbb{D})$, then $$\tilde{\Delta}\big(B_{\alpha}u\big) = 4(\alpha+1)(\alpha+2)(B_{\alpha}u-B_{\alpha+1}u).$$
This shows that for such $u$, $B_{\alpha+1}u = \Big(1-\dfrac{\tilde{\Delta}}{4(\alpha+1)(\alpha+2)}\Big)B_{\alpha}u$.
\end{lemma}

Lemma \ref{L:weightedBerezin} implies that for any integer $k\geq 1$, we have $B_{k}u = q_{k}(\tilde{\Delta})(B_{0}u)$ for $u\in L^1(\mathbb{D})$, where $q_{k}(\lambda) = \prod_{j=1}^{k}\Big(1-\frac{\lambda}{4j(j+1)}\Big).$

It was showed by Ahern and {\v C}u{\v c}kovi{\'c} in \cite{AhernJFA2001} that the Berezin transform $B_0$ and the invariant Laplacian $\tilde{\Delta}$ commute in the following sense: if $u$ is in $C^2(\mathbb{D})$ such that both $u$ and $\tilde{\Delta}u$ belong to $L^1(\mathbb{D})$, then $\tilde{\Delta}B_{0}u = B_{0}\tilde{\Delta}u$. By applying this identity multiple times, we see that if $u$ is $2s$ times continuously differentiable such that $u$ and $(\tilde{\Delta})^{j}u$ belong to $L^1(\mathbb{D})$ for all $1\leq j\leq s$, then
\begin{equation}\label{Eqn:commutingDeltaB_0}
q(\tilde{\Delta})B_0u = B_0\big(q(\tilde{\Delta})u\big)
\end{equation}
for any polynomial $q$ of degree $s$.

In the study of finite sums of finite products of Toeplitz operators with bounded harmonic symbols, one often encounters the problem of classifying integrable functions $u$ for which $B_0(u)$ satisfies certain conditions. The following theorem due to N.V. Rao \cite{Rao2010} shows that for locally bounded $u$, $B_0u$ cannot be of the form $f_1\bar{g}_1+\cdots +f_m\bar{g}_m$ with analytic functions $f_1,\ldots, f_m$ and $g_1,\ldots, g_m$ unless $u$ is harmonic. The proof uses recent results on finite rank Toeplitz operators (with distributional symbols) on the analytic Bergman space. In fact Rao considered $u\in L^1(\mathbb{D})$ without any boundedness assumption but the conclusion is different and we do not need it here for our goal.

\begin{theorem}\label{T:finiteRankBer} Suppose $u$ is an integrable function that is locally bounded on $\mathbb{D}$ such that $B_0(u)=f_1\bar{g}_1+\cdots +f_m\bar{g}_m$ for some positive integer $m$, where $f_1,\ldots, f_m$ and $g_1,\ldots, g_m$ are analytic on $\mathbb{D}$. Then $u$ is harmonic, or equivalently, $\tilde{\Delta}u=0$ on $\mathbb{D}$, and we have $u=f_1\bar{g}_1+\cdots+f_m\bar{g}_m$.
\end{theorem}

\begin{remark}
For analytic functions $f_1,\ldots, f_m$ and $g_1,\ldots, g_m$, Theorem 3.3 in \cite{ChoeRMI2008} shows that the function $f_1\bar{g}_1+\cdots+f_m\bar{g}_m$ is harmonic if and only if $(f_1-f_1(0))(\bar{g}_1-\bar{g}_1(0))+\cdots+(f_m-f(0))(\bar{g}_m-\bar{g}_m(0))=0$ on $\mathbb{D}$. In particular, when $m=1$, we have $(f_1-f_1(0))(g_1-g_1(0))=0$, which implies either $f_1$ or $g_1$ is a constant function.
\end{remark}


\section{Toeplitz operators on $A^2_n$}\label{S:ToeplitzOp}

Recall that for any function $f\in L^1(\mathbb{D},\mathrm{dA})$, the Toeplitz operator $T_f$ is defined by the formula
\begin{equation}\label{Eqn:defToeplitz}
(T_fh)(z) = \int_{\mathbb D}f(w)h(w)\bar{K}_z(w)\mathrm{dA}(w),
\end{equation}
for $h\in A^2_n$ for which the integral on the right hand side is well defined for all $z\in\mathbb{D}$. Since the function $\bar{K}_z(w)$ is $n$-analytic in the variable $z$, $T_fh$ is $n$-analytic on $\mathbb{D}$. Note that the domain of $T_f$ contains the space $\mathcal{W}$ of all linear combinations of kernel functions $K_z$, which is dense in $A^2_n$. We say that $T_f$ is bounded on $A^2_n$ if there is a constant $C>0$ such that for any $h\in\mathcal{W}$, $T_fh$ belongs to $A^2_n$ and $\|T_fh\|\leq C\|h\|$. In this case $T_f$ extends uniquely to a bounded operator on $A^2_n$. Also for $h\in A^2_n$ and $g=\sum_{j=1}^{m}\alpha_j K_{z_j}\in\mathcal{W}$, we have
\begin{align}\label{Eqn:innerProduct}
\langle T_fh,g\rangle & = \sum_{j=1}^{m}\bar{\alpha}_j\langle T_fh, K_{z_j}\rangle = \sum_{j=1}^{m}\bar{\alpha}_j (T_fh)(z_j)\notag\\
& = \sum_{j=1}^{m}\bar{\alpha}_j\int_{\mathbb{D}}f(w)h(w)\bar{K}_{z_j}(w)\mathrm{d}A(w)\quad \text{ (by \eqref{Eqn:defToeplitz})}\notag\\
& =\int_{\mathbb{D}}f(w)h(w)\bar{g}(w)\mathrm{d}A(w).
\end{align}

Let $P$ denote the orthogonal projection of $L^2(\mathbb{D})$ onto $A^2_n$. If $f$ belongs to $L^2(\mathbb{D})$, then for any bounded $h$ in $A^2_n$, the function $fh$ belongs to $L^2(\mathbb{D})$. This implies that $$P(fh)(z) = \langle P(fh),K_z\rangle =\langle fh, K_z\rangle = \int_{\mathbb{D}}f(w)h(w)\bar{K}_z(w)\mathrm{dA}(w)=(T_fh)(z).$$ We then recover the usual definition of $T_f$, which is given by $T_fh = P(fh)$.

If $f$ is a bounded function, then $T_f$ is a bounded operator on $A^2_n$ and we have $\|T_f\|\leq\|f\|_{\infty}$ and $T^{*}_{f}=T_{\bar{f}}$.

If $f$ is analytic (not necessarily bounded), then $T_f$ is just the operator of multiplication by $f$. Also for any $z, w$ in $\mathbb{D}$, we have
\begin{align*}
\langle T_{\bar{f}}K_z,K_w\rangle & = \langle K_z,fK_w\rangle = \bar{f}(z)\bar{K}_w(z)= \langle \bar{f}(z)K_z, K_w\rangle.
\end{align*}
Using the density in $A^2_n$ of the linear span of $\{K_w: w\in\mathbb{D}\}$, we conclude that $T_{\bar{f}}K_z=\bar{f}(z)K_z$ for $z\in\mathbb{D}$.

Suppose $u$ is a harmonic function that belongs to $L^2(\mathbb{D})$. Then there are analytic functions $f, g$ in $L^2(\mathbb{D})$ such that $u=f+\bar{g}$ and hence,
\begin{align*}
\langle T_uK_z, K_w\rangle & = \langle \bar{g}K_z,K_w\rangle+\langle K_z,\bar{f}K_w\rangle = (f(w)+\bar{g}(z))\langle K_z,K_w\rangle
\end{align*}
for all $z, w$ in $\mathbb{D}$. This shows, in particular, that $\langle T_uK_z, K_z\rangle = u(z)\|K_z\|^2$ for any $z\in\mathbb{D}$. Hence, if $T_u$ is a bounded operator, then $\|u\|_{\infty}\leq \|T_u\|$. On the other hand, as we mentioned above, $\|T_u\|\leq\|u\|_{\infty}$ whenever $u$ is bounded. We then conclude that $\|T_u\|=\|u\|_{\infty}$.

Using the kernel function $K(z,w)$, we define the Berezin transform of a bounded operator $T$ on $A^2_n$ as $B(T)(z) = \|K_z\|^{-2}\langle TK_z,K_z\rangle = \langle Tk_z,k_z\rangle$ for $z\in \mathbb{D}$. (Recall that $k_z=\|K_z\|^{-1}K_z$ is the normalized reproducing kernel at $z$.) It follows from the formula for $K_z$ that $B$ maps the space $\mathcal{B}(A^2_n)$ of bounded linear operators on $A^2_n$ to the space of real analytic functions on $\mathbb{D}$. Also $\|B(T)\|\leq \|T\|$ for all $T\in\mathcal{B}(A^2_n)$.

Suppose $f\in L^1(\mathbb{D})$ such that $T_f$ is a bounded operator on $A^2_{n}(\mathbb D)$. Then for $z\in\mathbb{D}$, using the formula for $k_z$ as in \eqref{Eqn:normalizedKernel}, we have
\begin{align*}
B(T_f)(z) & = \langle T_f k_z,k_z\rangle = \int_{\mathbb{D}} f(w)|k_z(w)|^2\mathrm{d}A(w)\quad\text{ (by \eqref{Eqn:innerProduct})}\\
& = \int_{\mathbb{D}}f(w)\dfrac{(1-|z|^2)^2}{|1-w\bar{z}|^{4}}\Big\{\sum_{j=0}^{n-1}(-1)^{j}\binom{n}{j+1}\binom{n+j}{n}|\varphi_{z}(w)|^{2j}\Big\}^2\mathrm{d}A(w).
\end{align*}
Note that the last integral is defined for all $z$ in $\mathbb{D}$ and all $f$ in $L^1(\mathbb{D})$ regardless of the boundedness of $T_f$. We will denote it by $Bf(z)$ and call the function $Bf$ the Berezin transform of $f$.

\subsection{Compactness of Toeplitz operators with non-negative symbols}

Compactness and Schatten class membership of Toeplitz operators (on analytic Bergman spaces) with non-negative symbols can be characterized using the Berezin transform of the symbol. See \cite[Chapter 7]{ZhuAMS2007} for more details. We will show that this characterization (for simplicity, we will only work on the compactness characterization) carries over to $n$-analytic Bergman spaces with a similar approach.

\begin{lemma}\label{L:inequalityBer}
There is a constant $c>0$ so that for any non-negative function $f\in L^1(\mathbb{D})$ such that $\tilde{f}=Bf$ belongs to $L^1(\mathbb{D})$ and $T_{\tilde{f}}$ is bounded on $A^2_n$, we have $T_{f}\leq cT_{\tilde f}$.
\end{lemma}

\begin{proof}
Let $h$ be a linear combination of kernel functions. We have
{\allowdisplaybreaks
\begin{align*}
\langle T_{\tilde{f}}h,h\rangle & = \int_{\mathbb{D}}\tilde{f}(z)|h(z)|^2\mathrm{dA}(z)\quad \text{ (by \eqref{Eqn:innerProduct})}\\
& = \int_{\mathbb{D}}\int_{\mathbb{D}}f(w)|k_z(w)|^2|h(z)|^2\mathrm{dA}(w)\mathrm{dA}(z)\\
& = \int_{\mathbb{D}}\Big(\int_{\mathbb{D}}|k_z(w)|^2|h(z)|^2\mathrm{dA}(z)\Big)f(w)\mathrm{dA}(w)\\
& \geq\frac{1}{c}\int_{\mathbb{D}}|h(w)|^2f(w)\mathrm{dA}(w)\quad \text{ (by Lemma \ref{L:integralEstimate})}\\
& = \frac{1}{c}\langle T_fh, h\rangle.
\end{align*}}
Since the above inequalities hold for all $h$ in a dense subset of $A^2_n$ and $T_{\tilde f}$ is bounded, we conclude that $T_f\leq cT_{\tilde{f}}$.
\end{proof}


Suppose $f$ is a bounded function whose support is a compact subset of the unit disk. Since $T_f$ is an integral operator with kernel $f(w)\bar{K}_{z}(w)$, which is bounded in $(z,w)\in\mathbb{D}\times\mathbb{D}$, we see that $T_f$ is a Hilbert-Schmidt operator, hence compact. Using an approximation argument, we obtain the next lemma, which is well known in the case of analytic Bergman space $A^2_{1}$.

\begin{lemma}\label{L:cptToeplitz}
Suppose $f$ is a bounded function such that $\lim\limits_{|z|\uparrow 1}f(z)=0$. Then the operator $T_f$ is compact on $A^2_{n}$.
\end{lemma}

We are now able to characterize compact Toeplitz operators on $A^2_n$ with non-negative symbols. These results are analogous to the corresponding results on the analytic Bergman space $A^2_{1}$.

\begin{theorem}\label{T:cptToeplitz} Let $f\in L^1(\mathbb{D})$ be a non-negative function. Then $T_f$ is compact on $A^2_n$ if and only if $Bf(z)\to 0$ as $|z|\uparrow 1$.
\end{theorem}

\begin{proof} 
Suppose $T_f$ is compact. Since $k_z\to 0$ weakly as $|z|\uparrow 1$ by Lemma \ref{L:behaviorPoly}, we have $Bf(z)=B(T_f)(z)=\langle T_fk_z,k_z\rangle\to 0$ as $|z|\uparrow 1$. 

Conversely, suppose $\tilde{f}(z)=Bf(z)\to 0$ as $|z|\uparrow 1$. This, in particular, implies that $\tilde{f}$ is bounded. Lemma \ref{L:cptToeplitz} then shows that $T_{\tilde{f}}$ is compact. Using Lemma \ref{L:inequalityBer}, we conclude that $T_{f}$ is also compact.
\end{proof}

\begin{remark}
It is possible to introduce the notion of Carleson measure on $A^2_n$ and use it to obtain Theorem \ref{T:cptToeplitz} as it was done for the analytic Bergman space $A^2_1$ in \cite[Section 7.2]{ZhuAMS2007}. We leave the details to the interested reader.
\end{remark}

\begin{remark}
On the analytic Bergman space $A^2_{1}(\mathbb D)$, Axler and Zheng showed (see the equivalence of (a) and (c) in \cite[Theorem 7.23]{ZhuAMS2007}) that for a bounded function $f$, the Toeplitz operator $T_f$ is compact if and only if $B_0(f)(z)\to 0$ as $|z|\uparrow 1$. We do not know if this result can be extended to any $A^2_n$ with $n\geq 2$. (Axler-Zheng's proof, when applied to $A^2_n$ with $n\geq 2$, faces many difficulties which cannot be resolved.) Note that in Theorem \ref{T:cptToeplitz}, we do not require the boundedness of $f$ but we do require the positivity, which is crucial for the proof to work.
\end{remark}

\subsection{Toeplitz operators with harmonic symbols}

Suppose $u$ and $v$ are two bounded harmonic functions such that $u=f+\bar{g}$ and $v=h+\bar{k}$, where $f, g, h, k$ are analytic. For any $z, w$ in $\mathbb{D}$, we have
\begin{align}\label{Eqn:harmonicToeplitz}
\langle T_uT_v K_z, K_w\rangle & = \langle T_vK_z, T_{\bar{u}}K_w\rangle\notag\\
& = \langle hK_z+\bar{k}(z)K_z, \bar{f}(w)K_w+gK_w\rangle\\
& = \big(f(w)h(w)+f(w)\bar{k}(z)+\bar{g}(z)\bar{k}(z)\big)\langle K_z,K_w\rangle + \langle \bar{g}hK_z,K_w\rangle.\notag 
\end{align}
In particular, by letting $w=z$ and dividing by $\|K_z\|^{2}$ we obtain
\begin{align}\label{Eqn:BerProduct}
B(T_uT_v) & = fh+f\bar{k}+\bar{g}\bar{k}+B(T_{\bar{g}h}) = uv-\bar{g}h+B(T_{\bar{g}h}).
\end{align}
For two non-zero functions $x$ and $y$ in $A^2_n$, the rank-one operator $x\otimes y$ is defined by $(x\otimes y) (f)=\langle f, y\rangle x$ for $f\in A^2_n$. For $z, w$ in $\mathbb{D}$, we have
\begin{align*}
\langle (x\otimes y)K_z, K_w\rangle & = \langle K_z, y\rangle\langle x, K_w\rangle = x(w)\bar{y}(z).
\end{align*}
This implies
\begin{align}\label{Eqn:BerRankOne}
B(x\otimes y)(z) = \frac{(1-|z|^2)^{2}}{n^2}x(z)\bar{y}(z).
\end{align}

The Berezin transform is injective on the space of bounded linear operators on $A^2_1$: for any $\alpha\geq 0$, if $B_{\alpha}(T)$ is the zero function on $\mathbb{D}$, then $T=0$ (see \cite[Proposition 6.2]{ZhuAMS2007}). It turns out that this is not the case on $A^2_{n}$ when $n\geq 2$. In fact, if $f(z)=z$, then $f$ and $\bar{f}$ belong to $A^2_n$ and by \eqref{Eqn:BerRankOne}, $B(f\otimes 1)=B(1\otimes\bar{f})$ but it is clear that $f\otimes 1\neq 1\otimes\bar{f}$. We mention in passing that this example also serves as an example to show that the Berezin transform is not injective on bounded linear operators acting on the harmonic Bergman space of the unit disk. Examples involving Toeplitz operators on $A^2_n$ can also be produced as follows.

\begin{example}[Non-injectivity of Berezin transform] Choose analytic monomials $u_1, u_2, u_3$ and conjugate analytic monomials $v_1, v_2, v_3$ such that $$\dfrac{(1-\bar{z}w)^{2}}{n^2} = u_1(w)v_1(z)+u_{2}(w)v_{2}(z)+u_3(w)v_3(z).$$
Using \eqref{Eqn:BerProduct} and \eqref{Eqn:BerRankOne} together with the linearity of the Berezin transform, we obtain
$$B(\sum_{j=1}^{3}T_{u_j}T_{v_j})(z) = \sum_{j=1}^{3}u_j(z)v_j(z) = \dfrac{(1-|z|^2)^2}{n^2} = B(1\otimes 1)(z)$$
for any $z\in\mathbb{D}$. On the other hand, from \eqref{Eqn:harmonicToeplitz},
\begin{align*}
\langle\Big(\sum_{j=1}^{3}T_{u_j}T_{v_j}\Big)K_z,K_w\rangle & = \Big(\sum_{j=1}^{3}u_j(w)v_j(z)\Big)\langle K_z,K_w\rangle = \dfrac{(1-\bar{z}w)^2}{n^2}K_z(w),\\
\langle (1\otimes 1)K_z,K_w\rangle & = \langle K_z,1\rangle\langle 1,K_w\rangle = 1,
\end{align*}
for all $z$ and $w$ in $\mathbb{D}$. Since $K_z(w)$ is not the same as $\frac{n^2}{(1-\bar{z}w)^2}$ for $n\geq 2$, we conclude that $T_{u_1}T_{v_1}+T_{u_2}T_{v_2}+T_{u_3}T_{v_3}\neq 1\otimes 1$.
\end{example}

Even though the Berezin transform is not injective on the space of bounded operators on $A^2_n$, the implication $Bf=0\Longrightarrow f=0$ is in fact valid if appropriate conditions are imposed on $f$ (in particular, it is true for all bounded functions $f$). We will see this by analyzing the relation between $B$ and $B_0$, the unweighted Berezin transform. Recall that for $f\in L^1$ and $z\in\mathbb{D}$, the formula for $Bf(z)$ has the form
\begin{align*}
Bf(z) & = \int_{\mathbb{D}}f(w)\dfrac{(1-|z|^2)^2}{|1-w\bar{z}|^{4}}\Big\{\sum_{j=0}^{n-1}(-1)^{j}\binom{n}{j+1}\binom{n+j}{n}|\varphi_{z}(w)|^{2j}\Big\}^2\mathrm{d}A(w)\\
& = \int_{\mathbb{D}}f(\varphi_{z}(\zeta))\Big\{\sum_{j=0}^{n-1}(-1)^{j}\binom{n}{j+1}\binom{n+j}{n}|\zeta|^{2j}\Big\}^2\mathrm{d}A(\zeta),
\end{align*}
by the change of variable $w=\varphi_z(\zeta)$. Let $\mu$ be the polynomial of degree $2n-2$ defined by $$\mu(t) = \Big\{\sum_{j=0}^{n-1}(-1)^j\binom{n}{j+1}\binom{n+j}{n}t^{j}\Big\}^2.$$
Rewriting $\mu(t)=b_0+b_1(1-t)+\cdots+b_{2n-2}(1-t)^{2n-2}$ and using the remark after Lemma \ref{L:weightedBerezin} we see that $Bf$ can be written as

\begin{align}\label{Eqn:weightedBer}
Bf & = b_0B_{0}(f)+\frac{b_1}{2}B_1(f)+\cdots+\frac{b_{2n-2}}{2n-1}B_{2n-2}(f)\notag\\
& = b_0B_0(f)+\frac{b_1}{2}q_1(\tilde{\Delta})B_0(f)+\cdots+\frac{b_{2n-2}}{2n-1}q_{2n-2}(\tilde{\Delta})B_0(f)\notag\\
& = Q(\tilde{\Delta})B_0(f),
\end{align}
where $Q(t)=b_0+\frac{b_1}{2}q_1(t)+\cdots+\frac{b_{2n-2}}{2n-1}q_{2n-2}(t)$ and $q_k(t)=\prod_{j=1}^{k}(1-\frac{t}{4j(j+1)})$ for $1\leq k\leq 2n-2$. The formula $Bf = Q(\tilde{\Delta})B_0f$ implies that in order to understand $B$, we need to understand $Q(\tilde{\Delta})$ and $B_0$.

We begin with a result regarding the eigenspaces of the invariant Laplacian. For any complex number $\lambda$, let $X_{\lambda}$ denote the linear space of all twice differentiable functions $u$ on $\mathbb{D}$ such that $\tilde{\Delta}u=\lambda u$. Theorem 4.2.7 in \cite{RudinSpringer1980} shows that $X_{\lambda}\cap L^{\infty}(\mathbb{D})\neq\emptyset$ if and only if $\lambda$ belongs to the set
\begin{align}\label{Eqn:eigenvalues}
\Omega_{\infty} & = \{\lambda\in\mathbb{C}: 4{\rm Re}\ \lambda+({\rm Im}\ \lambda)^2\leq 0\}.
\end{align}
This shows that if $u\in C^2(\mathbb{D})\cap L^{\infty}(\mathbb{D}), \lambda\notin\Omega_{\infty}$ and $\tilde{\Delta}u=\lambda u$, then $u=0$. As a consequence, we obtain
\begin{lemma}\label{L:injectiveInvLap}
Let $q$ be a polynomial of degree $s\geq 0$ whose roots lie outside $\Omega_{\infty}$. Suppose $u$ is a function which is $2s+2$ times continuously differentiable on $\mathbb{D}$ such that $q(\tilde{\Delta})\tilde{\Delta}u=0$ and $(\tilde{\Delta})^{j}u$ is bounded for all $1\leq j\leq s$. Then $\tilde{\Delta}u=0$, that is, $u$ is harmonic.
\end{lemma}

We now need to study the locations of the roots of the polynomial $Q(t)$. Because of the complexity of $Q$, we have not fully understood the locations of its roots for all $n\geq 2$. With the help of Maple, we are able to show that for $2\leq n\leq 25$, all roots of $Q$ lie outside $\Omega_{\infty}$. As a result, Lemma \ref{L:injectiveInvLap} can be applied to $Q$. We illustrate here the case $n=2$ in which the computation can be carried out by hand. With $n=2$, we have $\mu(t)=(2-3t)^2 = 1-6(1-t)+9(1-t)^2$ and hence $$Q(t)=1-3q_1(t)+3q_2(t) = 1-\frac{t}{8}+\frac{t^2}{64}.$$ Since the roots of $Q$ are $t=4\pm 4\sqrt{-3}$, they lie outside $\Omega_{\infty}$. Note that in this case all the roots have positive real parts, however this is not true for other values of $n$.

Motivated by the conditions in Lemma \ref{L:injectiveInvLap}, we define $\mathcal{S}$ to be the space of all smooth functions $g$ on $\mathbb{D}$ such that $z\mapsto (1-|z|^2)^{k+l}\frac{\partial^{k+l}g}{\partial z^{k}\partial\bar{z}^{l}}(z)$ is bounded on $\mathbb{D}$ for all $0\leq k+l$. Using Leibniz's rule, it can be verified that $\mathcal{S}$ is an algebra, that is, $g_1g_2$ belongs to $\mathcal{S}$ whenever $g_1$ and $g_2$ belong to $\mathcal{S}$. Leibniz's rule also implies that $\mathcal{S}$ is invariant under the action of $\tilde{\Delta}$, that is, $\tilde{\Delta}g$ belongs to $\mathcal{S}$ whenever $g$ belongs to $\mathcal{S}$.

Recall that the Bloch space $\mathcal{B}$ consists of all analytic functions $f$ on $\mathbb{D}$ such that the function $z\mapsto (1-|z|^2)f'(z)$ is bounded on $\mathbb{D}$. Theorem 5.4 in \cite{ZhuAMS2007} shows that if $f$ belongs to $\mathcal{B}$ then $(1-|z|^2)^{m}f^{(m)}(z)$ is bounded on $\mathbb{D}$ for all $m\geq 1$. It follows, by induction, that functions of the form $(1-|z|^2)^{m}f^{(m)}(z)$ $(m\geq 1)$, where $f\in\mathcal{B}$, belong to $\mathcal{S}$.

If $u$ is a bounded harmonic function on $\mathbb{D}$ then as well known, there are functions $f$ and $g$ in $\mathcal{B}$ such that $u=f+\bar{g}$. It then follows that $u$ belongs to $\mathcal{S}$ and hence, $\mathcal{S}$ contains all finite sums of finite products of bounded harmonic functions.

\begin{remark}\label{R:injectiveInvLap}
Any function $u$ in $\mathcal{S}$ satisfies the hypothesis of Lemma \ref{L:injectiveInvLap}. Therefore if $q$ is a polynomial whose roots lie outside $\Omega_{\infty}$ and $q(\tilde{\Delta})\tilde{\Delta}u=0$, then $\tilde{\Delta}u=0$.
\end{remark}

We are now ready to prove our main theorem. As we discussed above, we need to restrict ourselves to $A^2_n$ with $1\leq n\leq 25$. We believe that this restriction is artificial but we are not able to remove it at this time.

\begin{theorem}\label{T:finitePerturbation}
Let $u_1,\ldots,u_d$ and $v_1,\ldots,v_d$ be bounded harmonic functions on $\mathbb{D}$. For each $j$ we write $u_j = f_j+\bar{g}_j$ and $v_j=h_j+\bar{k}_j$, where $f_j, g_j, h_j, k_j$ belong to the Bloch space. Let $F$ be a function in $\mathcal{S}$ (for example, $F$ may be taken to be a finite sum of finite products of bounded harmonic functions on $\mathbb{D}$). If 
\begin{equation}\label{Eqn:finiteRank}
\sum_{j=1}^{d}T_{u_j}T_{v_j}-T_F = \sum_{l=1}^{s}x_l\otimes y_l,
\end{equation} where $x_l, y_j$ belong to $A^2_n$, then the following statements hold
\begin{enumerate}[(i)]
\item $\sum_{j=1}^{d}h_j\bar{g}_j-F$ is harmonic on $\mathbb{D}$.
\item $\sum_{j=1}^{d}u_j(z)v_j(z) - F(z) = \frac{(1-|z|^2)^2}{n^2}\sum_{l=1}^{s}x_l(z)\bar{y}_l(z)$ for $z\in\mathbb{D}$.
\end{enumerate}
\end{theorem}

\begin{proof}
For any $w, z\in\mathbb{D}$, we have
\begin{align*}
\langle\Big(\sum_{l=1}^{s}x_l\otimes y_l\Big)k_z,k_w\rangle & = \sum_{l=1}^{s}\langle x_l, k_w\rangle\langle k_z,y_l\rangle\\
& = \sum_{l=1}^{s}x_l(w)\bar{y}_l(z)(\|K_z\|\|K_w\|)^{-1}.
\end{align*}
If \eqref{Eqn:finiteRank} holds, then by applying Berezin transform on both sides and using \eqref{Eqn:BerProduct} and the above identities, we obtain
\begin{align}\label{Eqn:Berezin}
& B(\sum_{j=1}^{d}h_j\bar{g}_j)(z)+\sum_{j=1}^{d}f_j(z)h_j(z)+f_j(z)\bar{k}_j(z)+\bar{g}_j(z)\bar{k}_j(z) - B(F)\notag\\
&\hspace{1in} = \|K_z\|^{-2}\sum_{l=1}^{s}x_l(z)\bar{y}_l(z) = \frac{(1-|z|^2)^2}{n^2}\sum_{l=1}^{s}x_l(z)\bar{y}_l(z)
\end{align}
for all $z\in\mathbb{D}$.

Put $u=\sum_{j=1}^{d}h_j\bar{g}_j-F$. Since $\tilde{\Delta}u$ belongs to $\mathcal{S}$, we may apply \eqref{Eqn:weightedBer} and then \eqref{Eqn:commutingDeltaB_0} to obtain $B(u)=Q(\tilde{\Delta})B_0u=B_0(Q(\tilde{\Delta})u)$. Since $x_l$ and $y_l$ are $n$-analytic, the function $n^{-2}(1-|z|^2)^2x_l(z)\bar{y}_{l}(z)$ may be written as a finite sum of products of the form $f\bar{g}$, where $f$ and $g$ are analytic on $\mathbb{D}$. Now, \eqref{Eqn:Berezin} implies that $B_0(Q(\tilde{\Delta})u)$ is also a finite sum of such products. Using Theorem \ref{T:finiteRankBer}, we conclude that $\tilde{\Delta}Q(\tilde{\Delta})u=0$, or equivalently, $Q(\tilde{\Delta})\big(\tilde{\Delta}u\big)=0$.

Since all roots of $Q$ lie outside $\Omega_{\infty}$ and $\tilde{\Delta}u$ belongs to $\mathcal{S}$, Remark \ref{R:injectiveInvLap} shows that $\tilde{\Delta}u=0$, which proves the first statement in the conclusion of the theorem. The second statement follows from \eqref{Eqn:Berezin} and the fact that harmonic functions are fixed by the Berezin transform.
\end{proof}

Theorem \ref{T:finitePerturbation} enjoys many interesting consequences which are analogous to the case of unweighted Bergman space $A^2_1$. The following results, for $A^2_1$, were obtained by Guo, Sun and Zheng in \cite[Section 4]{GuoIJM2007}. Using Theorem \ref{T:finitePerturbation}, we see that they also hold for $A^2_n$.

\begin{corollary}\label{C:finiteRankProduct}
If $u$ and $v$ are bounded harmonic functions on $\mathbb{D}$ such that $T_{u}T_{v}$ is of finite rank, then either $u$ or $v$ is the zero function.
\end{corollary}

\begin{proof}
Write $u=f+\bar{g}$ and $v=h+\bar{k}$, where $f, g, h, k$ belong to the Bloch space. Suppose $T_{u}T_{v} = \sum_{l=1}^{s}x_l\otimes y_l$ for some $s\geq 1$ and $x_l, y_l$ are in $A^2_n$. By Theorem \ref{T:finitePerturbation}, $h\bar{g}$ is harmonic and $$u(z)v(z) = \dfrac{(1-|z|^2)^{2}}{n^2}\sum_{l=1}^{s}x_l(z)\bar{y}_l(z)$$ for $z\in\mathbb{D}$. Since $\lim\limits_{|z|\uparrow 1}(1-|z|^2)\varphi(z)=0$ for all $\varphi$ belonging to $A^2_n$ by Lemma \ref{L:behaviorPoly}, we conclude that $\lim\limits_{|z|\uparrow 1}u(z)v(z)=0$.

Since $h\bar{g}$ is harmonic, either $h$ or $g$ must be a constant function. This implies that either $\bar{v}$ or $u$ is analytic on $\mathbb{D}$. It then follows that either $u$ or $v$ is the zero function (see the proof of \cite[Corollary 3.11]{ChoeRMI2008}).
\end{proof}

\begin{corollary}\label{C:finiteRankCommutator}
If $u$ and $v$ are non-constant bounded harmonic functions on $\mathbb{D}$ such that $[T_u, T_v]=T_{u}T_{v}-T_{v}T_{u}$ is of finite rank, then either both $u$ and $v$ are analytic, or both $\bar{u}$ and $\bar{v}$ are analytic, or $u=\alpha v+\beta$ for some complex numbers $\alpha$ and $\beta$.
\end{corollary}

\begin{proof} Write $u=f+\bar{g}$ and $v=h+\bar{k}$ as in the proof of Corollary \ref{C:finiteRankProduct}. By Theorem \ref{T:finitePerturbation} part $(i)$, $h\bar{g}-f\bar{k}$ is harmonic. The conclusion of the corollary follows from the same argument as the proof of Theorem 1 in \cite[p.10]{AxlerIEOT1991}.
\end{proof}

\begin{corollary}\label{C:finiteRankSemi}
If $u$ and $v$ are bounded harmonic functions on $\mathbb{D}$ such that the semi-commutator $[T_u,T_v)=T_uT_v-T_{uv}$ is of finite rank, then either $\bar{u}$ or $v$ is analytic.
\end{corollary}

\begin{proof}
Write $u=f+\bar{g}$ and $v=h+\bar{k}$ as in the proof of Corollary \ref{C:finiteRankProduct}. By Theorem \ref{T:finitePerturbation} part $(i)$, $h\bar{g}-uv$ is harmonic. This shows that $f\bar{k}$ is harmonic, which implies either $f$ or $k$ is a constant function. The conclusion of the corollary then follows.
\end{proof}

\begin{corollary}\label{C:productToeplitz}
If $u$ and $v$ are bounded harmonic functions and $F$ is a function in $\mathcal{S}$ such that $T_uT_v=T_F$, then either $\bar{u}$ or $v$ is analytic and $F=\bar{u}v$.
\end{corollary}

\begin{proof}
As before, write $u=f+\bar{g}$ and $v=h+\bar{k}$. By Theorem \ref{T:finitePerturbation}, $F=uv$ and $h\bar{g}-F$ is harmonic, which implies that $f\bar{k}$ is harmonic. Therefore either $f$ or $k$ is a constant function. The conclusion of the corollary then follows.
\end{proof}

\begin{remark}
In \cite{AhernJFA2004}, Ahern obtained a version of Corollary \ref{C:productToeplitz} when $n=1$ under a weaker assumption that $F$ is merely bounded. We do not know whether this still holds true for $n\geq 2$.
\end{remark}

\section{Open questions}
In this section we list some questions that we have not been able to answer. We hope these questions will stimulate more research in this area.

Theorem \ref{T:finitePerturbation} and its corollaries are restricted by the assumption that $n\leq 25$. This is due to the fact that we have not found a way to verify whether the polynomial $Q$ has all roots outside of $\Omega_{\infty}$ for all $n$.

\begin{question} Is it true that all roots of $Q$ lie outside $\Omega_{\infty}$ for all $n\geq 2$?
\end{question}

Theorem \ref{T:cptToeplitz} is restricted to Toeplitz operators with non-negative symbols.  The well known Axler-Zheng theorem states that if $T$ is a finite sum of finite products of Toeplitz operators with bounded symbols on $A^2_1$ such that $B(T)(z)\to 0$ as $|z|\uparrow 1$, then $T$ is compact. (Recall that $B(T)$ is the Berezin transform of $T$ as defined in Section \ref{S:ToeplitzOp}.) Since the Berezin transform $B$ is no longer one-to-one on the space of bounded operators on $A^2_n$ when $n\geq 2$, it seems that the full version of Axler-Zheng theorem may not hold. On the other hand, the answer to the following question may be affirmative.

\begin{question} Let $f$ be bounded on $\mathbb{D}$. If $B(T_f)(z)\to 0$ as $|z|\uparrow 1$, does it follow that $T_f$ is compact on $A^2_n$ for $n\geq 2$?
\end{question}

The general zero-product problem asks: if $f$ and $g$ are bounded functions such that $T_fT_g=0$, does it follow that one of the functions must be zero? Corollary \ref{C:finiteRankProduct} answers this question in the affirmative when both functions are harmonic. In \cite{AhernActSM2004}, Ahern and the first author  also obtained the affirmative answer for operators on $A^2_1$ when one of the function is radial. In \cite{LeIEOT2009,LeCAOT2010}, the second author generalized this result to more than two functions all of which, except possible one, are radial. While the general zero-product problem (even on $A^2_1$) is still far from its solution, an answer to the following question may be possible. (We note that Toeplitz operators with radial symbols are diagonal with respect to the standard orthonormal basis of $A^2_1$. This, however, is not the case on $A^2_n$ with $n\geq 2$.)

\begin{question} Let $f$ and $g$ be bounded functions, one of which is radial. If $T_fT_g=0$ on $A^2_n$ (or more generally, $T_fT_g$ has finite rank), must one of these functions be zero?
\end{question}

The next question is related to Corollary \ref{C:productToeplitz}. The answer was known to be affirmative on $A^2_1$ by Ahern \cite{AhernJFA2004}.

\begin{question} Suppose $u$ and $v$ are bounded harmonic functions and $T_uT_v=T_F$ on $A^2_n$ ($n\geq 2$), where $F$ is a bounded function (without any assumption about smoothness). Does it imply that $\bar{u}$ or $v$ is analytic?
\end{question}


\end{document}